\documentclass[12pt,twoside]{article}    
\usepackage{amsfonts,amsmath,latexsym}

\newtheorem{theorem}{Theorem}[section]

\newtheorem{lemma}[theorem]{Lemma}

\numberwithin{equation}{section}
\newenvironment{proof}[1][Proof]{\textbf{#1.} }{\ \rule{0.5em}{0.5em}}

\makeatletter \@addtoreset{equation}{section} \makeatother
\begin{document}


\pagestyle{myheadings}

\markboth{\hfill {\small AbdulRahman Al-Hussein } \hfill}{\hfill
{\small Sufficient conditions for optimality for SEEs} \hfill }


\thispagestyle{plain}


\begin{center}
{\large \textbf{Sufficient conditions for optimality for stochastic evolution equations}$^{*}$\footnotetext{$^{*}$
This work is supported by the Science College Research Center at Qassim University, project no. SR-D-012-1282. }} \\
\vspace{0.7cm} {\large AbdulRahman Al-Hussein }
\\
\vspace{0.2cm} {\footnotesize
{\it Department of Mathematics, College of Science, Qassim University, \\
 P.O.Box 6644, Buraydah 51452, Saudi Arabia \\ {\emph E-mail:} hsien@qu.edu.sa, alhusseinqu@hotmail.com}}
\end{center}

\begin{abstract}
In this paper we derive for a controlled stochastic
evolution system on a Hilbert space sufficient conditions for
optimality. Our result is derived by using its so-called adjoint backward
stochastic evolution equation.
\end{abstract}



\sloppy

\section{Introduction}\label{sec1}
Consider a stochastic controlled problem governed by the following stochastic evolution equation (SEE):
\begin{eqnarray*}
 \left\{ \begin{array}{ll}
              d X(t) = ( A X(t) + b (X(t), \nu (t)) ) d t
              + \sigma (X(t), \nu (t)) d W(t) , \;\; t \in [0 , T] ,\\
             \; \, X(0) = x_0 .
         \end{array}
 \right.
\end{eqnarray*}
This system is driven mainly by a possibly unbounded linear operator $A$ on a separable Hilbert space $H$ and a
cylindrical Wiener process $W$ on $ H .$ The control process here is denoted by $\nu (\cdot ) .$
Then the control problem is to minimize the cost functional, which is given by equation
(\ref{cost functional}) in Section~\ref{sec3}, over a set of admissible controls.

We shall concentrate in providing sufficient conditions for
optimality of this optimal control problem, which gives this minimization.
For this purpose we shall apply the theory of \emph{backward stochastic evolution
equations} (or shortly BSEE) as in equation (\ref{adjoint-bse}) in Section~\ref{sec3}, which together with backward stochastic differential
equations have become nowadays of great interests in many different fields. For example one can see
\cite{[H-Pe96]}, \cite{[Pe-93]}, \cite{[Oks05]}, \cite{[Oks-Zh05]} and \cite{[Y-Z]} for the
applications of backward stochastic differential equations in such optimal control problems.

Our work will not need studying Hamilton-Jacobi-Bellman equation either by using
semi-group technique or the technique of viscosity solutions. We refer the reader to \cite{[Cerrai_book]} and some of the related references therein for the semi-group technique.

Let us remark that necessary conditions for optimality of the control $\nu (\cdot )$ and its corresponding solution $X^{\nu (\cdot)}$ but for the case when $\sigma$ does not depend on $\nu$ can be found in \cite{[H-Pe96]}. This is also the case considered in our earlier work in \cite{[Alh-COSA]}. So the present paper generalize the work in \cite{[Alh-COSA]}.

\section{Notation}\label{sec2}
Let $(\Omega, \mathcal{F}, \mathbb{P})$ be a complete probability space and denote
by $\mathcal{N}$ the collection of $\mathbb{P}$\,-\,null sets of
$\mathcal{F} .$ Let $\{ W (t) ,\, 0 \leq t \leq T \} $ be a
cylindrical Wiener process on $H$ with its completed natural
filtration ${\mathcal{F}_t = \sigma \{ \ell \circ W (s) \, , \; 0
\leq s \leq t \, , \ell \in H^{*} \}\vee \mathcal{N} } , \; t \geq 0
;$ cf. \cite{[preprint_MRT]}.

For a separable Hilbert space $E$ let $L^2_{\mathcal{F}} (
0, T; E )$ denote the space of all $\{ \mathcal{F}_t , 0
\leq t \leq T \}$\,-\,progressively measurable processes $f$
with values in $E$ such that \[ \mathbb{E}\; [ \int_0^T |
f(t) |_{E}^2 \; dt ] < \infty .\] Thus
$L^2_{\mathcal{F}} ( 0, T; E )$ is a Hilbert space with the
norm \[ || f || \; = \Big{(} \mathbb{E}\; [ \int_0^T |
f(t) |_{E}^2 \; dt ] \Big{)}^{1/2}\, . \]

It is known as in \cite{[Da-Z]} that for $ f \in L^2_{\mathcal{F}} (
0, T; L_2 (H) ) ,$ where $L_2 (H)$ is the space of all
Hilbert-Schmidt operators on $H ,$  the stochastic integral $\int f(t)
dW(t)$ can be defined as a continuous stochastic process in $H .$ The inner product on $L_2 (H)$ will be denoted by by $\big{<} \cdot , \cdot \big{>}_2 .$

\section{Results}\label{sec3}
Let $\mathcal{O}$ be a separable Hilbert space equipped with an inner
product $\big{<} \cdot , \cdot \big{>}_{\mathcal{O}}$, and let $U$ be
a convex subset of $\mathcal{O} .$ We say that
$\nu (\cdot ) : [0 , T]\times \Omega \rightarrow \mathcal{O}$ is \emph{admissible}
if $\nu (\cdot ) \in  L^2_{\mathcal{F}} ( 0 , T ; \mathcal{O} )$ and $\nu (t) \in U \; \; a.e., \; a.s.$
The set of admissible controls will be denoted by $\mathcal{U}_{ad} .$
Let $b: H \times \mathcal{O} \rightarrow H$ and $\sigma: H \times \mathcal{O} \rightarrow L_2(H)$ be two continuous mappings. Consider the following controlled system:
\begin{eqnarray}\label{forward-see}
\left\{ \begin{array}{ll}
              d X(t) = ( A X(t) + b (X(t), \nu (t)) ) d t
              + \sigma (X(t), \nu (t))  d W(t) ,\\
             \; \, X(0) = x_0 ,
         \end{array}
 \right.
\end{eqnarray}
where $\nu (\cdot ) \in \mathcal{U}_{ad} $ represents a control variable.
A solution of $(\ref{forward-see})$ will be denoted by $X^{\nu (\cdot )}$ to indicate
the presence of the control.

Let $\ell : H \times \mathcal{O} \rightarrow \mathbb{R}$ and $\phi : H \rightarrow \mathbb{R}$ be two measurable mappings such that the following {\it cost functional} is defined:
\begin{equation}\label{cost functional}
J(\nu (\cdot ) ) : = \mathbb{E} \; [ \; \int_0^T \ell ( X^{\nu (\cdot )}
(t) , \nu (t) )  dt + \phi ( X^{\nu (\cdot )} (T) ) \; ] , \;\; \nu (\cdot ) \in \mathcal{U}_{ad} .
\end{equation}
For example we can take $\ell$ and $\phi$ to satisfy the assumptions given in Theorem~\ref{main thm}.

The optimal control problem of the system $(\ref{forward-see})$ is to find the {\it
value function } $J^*$ and \emph{an optimal control} $\nu^{*} (\cdot
) \in \mathcal{U}_{ad} $ such that
\begin{equation}\label{control problem}
J^{*} : = \inf \{ J(\nu (\cdot ) ) : \; \nu (\cdot ) \in \mathcal{U}_{ad} \} =
J(\nu^{*} (\cdot ) )  .
\end{equation}
If this happens, the corresponding solution $X^{\nu^{*} (\cdot )}$ is
called \emph{an optimal solution} of the stochastic control problem
(\ref{forward-see})--(\ref{control problem}) and
$( X^{\nu^{*} (\cdot )}\, , \nu^{*} (\cdot ) )$ is called \emph{an optimal pair.}

Let us now recall the following theorem.
\begin{theorem}\label{thm1}
Assume that $A$ is an unbounded linear operator on $H$ that
generates a $C_0$-semigroup $\{ S(t) ,\; t \geq 0 \} $ on $H$, and
$b , \sigma$ are continuously Fr\'echet differentiable with respect to
$x$ and their derivatives $b_x \, , \,  \sigma_x$ are uniformly bounded.
Then for every $\nu (\cdot ) \in \mathcal{U}_{ad} $ there exists a unique mild solution $X^{\nu (\cdot )}$
on $[0, T] $ to $(\ref{forward-see}) .$ That is $X^{\nu (\cdot )}$ is a progressively measurable stochastic
process such that $X (0) = x_0 $ and for all $t \in [0 , T ] ,$
\begin{eqnarray}\label{thm1:solution}
X^{\nu (\cdot )} (t) &=& S ( t ) x_0 + \int_0^t S ( t - s) b(X^{\nu (\cdot )} (s) , \nu (s) ) ds \nonumber  \\
& & \hspace{2cm} + \, \int_0^t S ( t - s) \, \sigma (X^{\nu (\cdot )} (s) , \nu (s) ) \, d W(s) .
\end{eqnarray}
\end{theorem}
The proof of this theorem can be derived in a similar way to those
in \cite[Chapter 7]{[Da-Za]} or \cite{[Ichi]}.

\smallskip

From here on we shall assume that $A$ is the generator of a $C_0$-semigroup $\{ S(t) ,\; t \geq 0 \} $ on $H .$
Its adjoint operator $A^{*} : \mathcal{D} ( A^{*} ) \subset H \rightarrow H$ is then the infinitesimal generator of
the adjoint semigroup $\{ S^{*} (t) ,\; t \geq 0 \} $ of $\{ S (t) \, , t \geq 0 \} . $

\smallskip

As it is known that backward stochastic differential equations play
an important role in deriving the maximum (or minimum) principle for SDEs, it is natural to search for such a role for SEEs like
(\ref{forward-see}). For this purpose, let us first consider the \emph{Hamiltonian}: $$ \mathcal{H} : H\times \mathcal{O} \times H
\times L_2(H) \rightarrow \mathbb{R} ,$$
\begin{eqnarray}\label{def:Hamiltonian}
\mathcal{H} (x, \nu , y , z ): = \ell (x , \nu) \, + \big{<} b(x , \nu ) , y\big{>}_{H} +
\big{<}\sigma (x , \nu ) , z \big{>}_2 .
\end{eqnarray}

Consider the following adjoint BSEE on $H$:
\begin{eqnarray}\label{adjoint-bse}
\left\{ \begin{array}{ll}
             -\, d Y^{\nu (\cdot )} (t) = & \big(\, A^{*} \, Y^{\nu (\cdot )} (t) + \nabla_{x} \mathcal{H}
             ( X^{\nu (\cdot )}(t), \nu (t), Y^{\nu (\cdot )}(t) , Z^{\nu (\cdot )} (t) ) \,
             \big)\, dt \\& \hspace{1.70in} - Z^{\nu (\cdot )} (t) d W(t) , \; \; 0 \leq t < T,  \\
             \; \; \; Y^{\nu (\cdot )} (T) = & \nabla \phi (X^{\nu (\cdot )}(T)),
         \end{array}
 \right.
\end{eqnarray}
where $\nabla \phi$ denotes the gradient of $\phi ,$ which is
defined, by using the directional derivative $D\phi (x) (h)$ of
$\phi$ at a point $x \in H$ in the direction of $h \in H ,$ as
$\big{<} \nabla \phi (x) , h \big{>}_{H} = D\phi (x) (h) \; ( \, = \phi_x (h) \, ).$

\bigskip

A \emph{mild solution} (or briefly a solution) of (\ref{adjoint-bse}) is a pair $( Y , Z ) \in L^2_{\mathcal{F}}( 0,
T; H ) \times L^2_{\mathcal{F}}( 0, T; L_2(H)) $ such that we have
$\mathbb{P}$\,-\,a.s. for all $t \in [0 , T]$
\begin{eqnarray}\label{bse:sol}
Y^{\nu (\cdot )} (t) &=& S^{*} ( T - t )\,  \nabla \phi (X^{\nu (\cdot )}(T)) \nonumber \\ &
& + \, \int_t^T S^{*} ( s - t )\, \nabla_{x} \mathcal{H} (
X^{\nu (\cdot )}(s), \nu (s), Y^{\nu (\cdot )}(s) , Z^{\nu (\cdot )} (s) ) ds \nonumber \\
& & \hspace{1.5in} - \, \int_t^T S^{*} ( s - t)\, Z^{\nu (\cdot )}(s) dW(s) .
\end{eqnarray}

Existence of such solutions can be derived from the following theorem.

\begin{theorem}[\cite{[BSEEs]} or \cite{[H-Pe91]}] \label{th:solution of adjointeqn}
Assume that $b , \sigma , \ell , \phi$ are continuously Fr\'echet differentiable with respect to
$x ,$ the derivatives $b_x , \sigma_x , \sigma_{\nu} , \ell_x $ are uniformly bounded, and
$$ | \phi_x |_{L(H,H)} \leq C \, ( 1 + |x|_{H} ) $$ for some constant $C > 0 .$

Then there exists a unique mild solution $( Y^{\nu (\cdot )} , Z^{\nu (\cdot )})$ of the
BSEE~(\ref{adjoint-bse}).
\end{theorem}

An alternative proof of this theorem by using finite dimensional framework through the Yosida approximation of $A$ can be found in \cite{[Tess96]}.

\bigskip

Now we state our main result.
\begin{theorem}\label{main thm}
For a given admissible control $\nu^* (\cdot )$ let $X^{\nu^{*} (\cdot )}$ and $( Y^{\nu^{*} (\cdot )} , Z^{\nu^{*} (\cdot )} )$ be the solutions of
the corresponding equations (\ref{forward-see}) and
(\ref{adjoint-bse}) respectively. Suppose that \\
(i)\; $\phi$ is convex, \\
(ii)\; $b , \sigma , \ell$ are continuously Fr\'echet differentiable with respect to
$x , \nu ,$ $\phi$ is continuously Fr\'echet differentiable with respect to
$x ,$ the derivatives $b_x , \, b_{\nu} , \, \sigma_x , \sigma_{\nu} , \ell_x , \, \ell_{\nu}$ are uniformly bounded, and
$$ | \phi_x |_{L(H,H)} \leq C \, ( 1 + |x|_{H} ) $$ for some constant $C > 0 ,$ \\
(iii)\; $\mathcal{H} ( \cdot , \cdot , Y^{\nu^{*} (\cdot )}(t) , Z^{\nu^{*}
}(t) )$ is convex for all $t \in [0 , T] , \; \mathbb{P}$\,-\,a.s., and \\
(iv)\; ${ \mathcal{H} ( X^{\nu^{*} (\cdot )} (t) , \nu^{*} (t) , Y^{\nu^{*} (\cdot )}
(t) , Z^{\nu^{*} (\cdot )} (t) ) = \displaystyle{ \inf_{\nu \in \, U } } \;
\mathcal{H} ( X^{\nu^{*} (\cdot )}(t) , \nu , Y^{\nu^{*} (\cdot )}(t) , Z^{\nu^{*}
}(t) ) }$ \\ for a.e. $t \in [0 , T] , \; \mathbb{P}$\,-\,a.s.

Then $( X^{\nu^{*} (\cdot )} \, , \nu^{*} (\cdot ) )$ is an optimal
pair for the problem (\ref{forward-see})--(\ref{control problem}) .
\end{theorem}


Examples stated in the introduction of \cite{[Tess96]} and \cite{[Alh-COSA]} are covered by this theorem.

\section{Proofs}\label{sec4}
In this section we shall establish the proof of Theorem~\ref{main thm}.
Keeping the notation in Theorem~\ref{main thm} we need the following duality relation between (\ref{forward-see}) and (\ref{adjoint-bse}).

\begin{lemma}\label{lem:lemma1}
Let $\psi_{1} (t) :=  b ( X^{\nu^{*} (\cdot )}(t) , \nu^{*} (t)) - b( X^{\nu (\cdot )}(t) , \nu (t)) ,$ where ${ t \in [0 , T] . }$ If assumption (ii) in Theorem~\ref{main thm} holds, then
\begin{eqnarray}\label{eq:eq5-2}
& & \mathbb{E} \; [\; \big{<} \; Y^{\nu^{*} (\cdot )} (T) \; ,
X^{\nu^{*} (\cdot )} (T) - X^{\nu (\cdot )} (T) \; \big{>} \; ] = \mathbb{E} \; [\; \int_0^T \big{<} \;
Y^{\nu^{*} (\cdot )} (t) \; , \psi_{1} (t) \; \big{>} \; dt \; ] \nonumber \\
& & - \; \mathbb{E} \; [\; \int_0^T  \big{<}\;
\nabla_{x} \mathcal{H} ( X^{\nu^{*} (\cdot )} (t), \nu^{*} (t), Y^{\nu^{*}
}(t) , Z^{\nu^{*} (\cdot )} (t) ) \, , X^{\nu^{*} (\cdot )} (t) - X^{\nu (\cdot )} (t) \;
\big{>} \; dt \; ] \nonumber
\\ &&
+ \;  \mathbb{E} \; [ \; \int_0^T \big< \sigma (
X^{\nu^{*} (\cdot )}(t) , \nu^{*} (t)) - \sigma ( X^{\nu (\cdot )}(t) , \nu (t)) \, , \, Z^{\nu^{*} (\cdot )} (t) \big>_2 \, dt \; ] .
\end{eqnarray}
\end{lemma}
\begin{proof}
The proof is done by using Yosida approximation of the operator $A$ and It\^{o}'s formula, and can be gleaned easily from the proof of Theorem~2.1 in \cite{[Tess96]}.
\end{proof}

\bigskip

We are now ready to establish the proof of Theorem~\ref{main thm}.

\bigskip

\noindent \begin{proof}[ {\bf Proof of Theorem~\ref{main thm}}]
Let $\nu (\cdot )$ be an arbitrary admissible control. From the
definitions in (\ref{control problem}) and (\ref{cost functional}) we
obtain
\begin{eqnarray}\label{pf-eq1}
J ( \nu^{*}  (\cdot ) ) - J ( \nu  (\cdot ) )
&=& \mathbb{E} \; [\; \int_0^T \ell ( X^{\nu^{*} (\cdot )}(t) , \nu^{*}
(t ) )  \, dt + \phi (  X^{\nu^{*} (\cdot )}(T) ) \; ] \nonumber \\
& & - \, \mathbb{E} \; [ \; \int_0^T \ell ( X^{\nu (\cdot )}(t) , \nu (t ) ) dt + \phi (  X^{\nu (\cdot )}(T) )  \; ] \nonumber \\
&=& \mathbb{E} \; [\; \int_0^T ( \ell ( X^{\nu^{*} (\cdot )}(t) ,  \nu^{*}
(t ) ) - \ell ( X^{\nu (\cdot )}(t) ,  \nu (t )  ) dt \; ] \nonumber \\
& & + \; \mathbb{E} \; [ \; \phi (  X^{\nu^{*} (\cdot )}(T) ) - \phi (
X^{\nu (\cdot )}(T) ) \; ] .
\end{eqnarray}
But
\begin{eqnarray*}
&& \hspace{-1cm} \ell ( X^{\nu^{*} (\cdot )}(t) ,  \nu^{*} (t ) ) - \ell ( X^{\nu (\cdot )}(t) ,  \nu (t ) ) =  \mathcal{H} (X^{\nu^{*} (\cdot )}(t),
\nu^{*} (t), Y^{\nu^{*} (\cdot )}(t) , Z^{\nu^{*} (\cdot )}(t) )
\\ & & \hspace{2.5cm} - \mathcal{H} (X^{\nu (\cdot )}(t), \nu (t), Y^{\nu^{*} (\cdot )}(t) , Z^{\nu^{*} (\cdot )} (t) ) \\
& & \hspace{2.5cm} - \, \big{<} b( X^{\nu^{*} (\cdot )}(t) , \nu^{*} (t) ) - b( X^{\nu (\cdot )}(t) , \nu (t) ) \, , \,
Y^{\nu^{*} (\cdot )} (t) \big{>} \\ & & \hspace{2.5cm} - \; \big{<} \sigma ( X^{\nu^{*} (\cdot )}(t) , \nu^{*} (t) ) - \sigma ( X^{\nu (\cdot )}(t) , \nu (t) ) \, , \, Z^{\nu^{*} (\cdot )} (t) \big{>}_2 \;\;\;\; a.s.
\end{eqnarray*}
Therefore (\ref{pf-eq1}) becomes
\begin{eqnarray}\label{pf-eq2}
&& \hspace{-0.9cm} J ( \nu^{*}  (\cdot ) ) - J ( \nu  (\cdot ) ) = \mathbb{E} \; \Big[\;
\int_0^T \Big( \, \mathcal{H} (X^{\nu^{*} (\cdot )}(t), \nu^{*} (t), Y^{\nu^{*}
}(t) , Z^{\nu^{*} (\cdot )}(t) ) \nonumber
\\ & &  \hspace{3cm} - \; \mathcal{H} (X^{\nu (\cdot )}(t), \nu (t), Y^{\nu^{*} (\cdot )}(t) , Z^{\nu^{*} (\cdot )} (t) ) \
\nonumber \\
& & \hspace{3cm} - \; \big{<} b ( X^{\nu^{*} (\cdot )}(t) , \nu^{*} (t) ) - b( X^{\nu (\cdot )}(t) , \nu (t) ) \; ,
Y^{\nu^{*} (\cdot )}(t) \big{>} \;
\nonumber \\
& &
\hspace{3cm} - \; \big{<} \sigma ( X^{\nu^{*} (\cdot )}(t) , \nu^{*} (t) ) - \sigma ( X^{\nu (\cdot )}(t) , \nu (t) ) \, , \,
Z^{\nu^{*} (\cdot )} (t) \big{>}_2
\Big)\, dt \; \Big]
\nonumber \\
& & \hspace{2.1in} + \; \mathbb{E} \; [ \; \phi (  X^{\nu^{*} (\cdot )}(T) )  - \phi (
X^{\nu (\cdot )}(T) ) \; ] .
\end{eqnarray}

By the convexity assumption on $\phi $ in (i) we get
\begin{eqnarray}\label{pf-eq3}
&& \hspace{-1.8cm} \phi (  X^{\nu^{*} (\cdot )}(T) )  - \phi (  X^{\nu (\cdot )}(T) ) \leq \; \big{<}
\; \nabla \phi (  X^{\nu^{*} (\cdot )}(T) ) \; , X^{\nu^{*} (\cdot )}(T) - X^{\nu (\cdot )}(T) \; \big{>} \; \; a.s.
\end{eqnarray}
But $\nabla \phi (  X^{\nu^{*} (\cdot )}(T) ) = Y^{\nu^{*} (\cdot )}(T) .$
Thus
\begin{eqnarray}\label{pf-eq4}
&& \hspace{-1.5cm} \mathbb{E} \, [\,  \phi (  X^{\nu^{*} (\cdot )}(T) )  - \phi (  X^{\nu (\cdot )}(T)
) \, ] \leq \; \mathbb{E} \, [ \, \big{<} \; Y^{\nu^{*} (\cdot )}(T) \; ,
X^{\nu^{*} (\cdot )} (T) - X^{\nu (\cdot )}(T) \; \big{>} \, ] .
\end{eqnarray}

Denote for $t \in [ 0 , T ],$
\begin{eqnarray*}
\delta \mathcal{H} (t) &=& \mathcal{H} (X^{\nu^{*} (\cdot )}(t), \nu^{*} (t),
Y^{\nu^{*} (\cdot )}(t) , Z^{\nu^{*} (\cdot )}(t) )
\\ & & - \mathcal{H} (X^{\nu (\cdot )}(t), \nu (t), Y^{\nu^{*} (\cdot )} (t) , Z^{\nu^{*} (\cdot )} (t)
)
\end{eqnarray*}
and
\begin{eqnarray*}
&& \psi_2 (t) = \delta \mathcal{H} (t) + \big{<} - \psi_1 (t) \; ,
Y^{\nu^{*} (\cdot )} (t) \, \big{>} \\ && \hspace{1in} - \; \big< \sigma (
X^{\nu^{*} (\cdot )}(t) , \nu^{*} (t)) - \sigma ( X^{\nu (\cdot )}(t) , \nu (t)) \, , \, Z^{\nu^{*} (\cdot )} (t) \big>_2 .
\end{eqnarray*}

By applying (\ref{pf-eq2}) and (\ref{pf-eq4}) we obtain
\begin{equation*}\label{pf-eq6}
J ( \nu^{*} (\cdot ) ) - J ( \nu (\cdot ) ) \leq
\mathbb{E} \; [\; \int_0^T \psi_2 (t) dt \; ] + \mathbb{E} \; [\; \big{<} \; Y^{\nu^{*} (\cdot )} (T) \; ,
X^{\nu^{*} (\cdot )} (T) - X^{\nu (\cdot )} (T) \; \big{>} \; ].
\end{equation*}
Consequently, by using Lemma~\ref{lem:lemma1} this inequality becomes
\begin{eqnarray}\label{pf-eq7}
& & \hspace{-1.25cm} J ( \nu^{*} (\cdot ) ) - J ( \nu (\cdot ) ) \leq \nonumber
\\ & &  - \, \mathbb{E} \; [\; \int_0^T
\big{<} \; \nabla_{x} \mathcal{H} ( X^{\nu^{*} (\cdot )} (t), \nu^{*} (t),
Y^{\nu^{*} (\cdot )}(t) , Z^{\nu^{*} (\cdot )} (t) ) \, , \, X^{\nu^{*} (\cdot )} (t) - X^{\nu (\cdot )}
(t) \; \big{>} \; dt \; ]  \nonumber \\ & & \hspace{-1.5cm}
  \hspace{3.75in} + \;
\mathbb{E} \; [\; \int_0^T \delta \mathcal{H} (t) \; dt \; ] .
\end{eqnarray}

From the convexity assumption of $\mathcal{H}$ in condition (iii) we see
that the following inequality holds for a.e. $t \in [ 0 , T ], \;
\mathbb{P}$\,-\,a.s.
\begin{eqnarray*}
\delta \mathcal{H} (t) &\leq&  \big{<} \, \nabla_{x} \mathcal{H}
(X^{\nu^{*} (\cdot )}(t), \nu^{*} (t), Y^{\nu^{*} (\cdot )}(t) , Z^{\nu^{*} (\cdot )}(t) ) \,
, \, X^{\nu^{*} (\cdot )} (t) - X^{\nu (\cdot )} (t) \, \big{>}
\\ & & + \; \big{<} \, \nabla_{\nu } \mathcal{H}
(X^{\nu^{*} (\cdot )}(t), \nu^{*} (t), Y^{\nu^{*} (\cdot )}(t) , Z^{\nu^{*} (\cdot )}(t) ) \,
, \, \nu^{*} (t) - \nu (t) \, \big{>}_{\mathcal{O}} .
\end{eqnarray*}
But the minimum condition (iv) implies
\begin{eqnarray*}
\big{<} \, \nabla_{\nu } \mathcal{H} (X^{\nu^{*} (\cdot )}(t), \nu^{*} (t),
Y^{\nu^{*} (\cdot )}(t) , Z^{\nu^{*} (\cdot )}(t) ) \, , \, \nu^{*} (t) - \nu (t) \,
\big{>}_{\mathcal{O}} \leq 0
\end{eqnarray*}
for a.e. $t \in [ 0 , T], \; \mathbb{P}$\,-\,a.s.; see e.g. \cite{[Lions]}.
Consequently, for a.e. $t \in [ 0 , T], \; \mathbb{P}$\,-\,a.s.,
\begin{eqnarray*}
\delta \mathcal{H} (t) \, - \big{<} \, \nabla_{x} \mathcal{H}
(X^{\nu^{*} (\cdot )}(t), \nu^{*} (t), Y^{\nu^{*} (\cdot )}(t) , Z^{\nu^{*} (\cdot )}(t) ) \,
, \, X^{\nu^{*} (\cdot )} (t) - X^{\nu (\cdot )} (t) \, \big{>} \; \leq 0 .
\end{eqnarray*}

Now by applying this result in (\ref{pf-eq7}) we deduce finally that $
{ J ( \nu^{*} (\cdot ) ) \leq J ( \nu (\cdot ) ). }$ This completes the
proof.
\end{proof}

\fussy


{\small

\begin{thebibliography}{99}
\bibitem{[preprint_MRT]} Al-Hussein, A.,  Martingale Representation Theorem In
Infinite Dimensions, \emph{Arab J. Math. Sc.}, 10, 1 (2004), 1--18.
\bibitem{[preprint-SW]} Al-Hussein, A.,  Strong, mild and weak solutions of
backward stochastic evolution equations, \emph{Random Oper. and
Stoch. Equ.,} 13, 2 (2005), 129--138.
\bibitem{[BSEEs]} Al-Hussein, A., Time-dependent Backward Stochastic Evolution
Equations, \emph{Bull. Malays. Math. Sci. Soc.,} 30, 2 (2007), 159--183.
\bibitem{[Alh-COSA]} Al-Hussein, A., Sufficient conditions of optimality for backward stochastic evolution equations,
Commun. Stoch. Anal. 4, 3 (2010), 433--442.
\bibitem{[Cerrai_book]} Cerrai, S., \emph{Second order PDE's in finite and infinite dimension.
A probabilistic approach.} Lecture Notes in Mathematics, 1762.
Springer-Verlag, Berlin, 2001.
\bibitem{[Da-Z]} Da Prato, G., Zabczyk, J., \emph{Stochastic equations
in infinite dimensions,} Encyclopedia of Mathematics and its
Applications, 44, Cambridge University Press, Cambridge, 1992.
\bibitem{[Da-Za]} Da Prato, G., Zabczyk, J., \emph{Second order partial differential
equations in Hilbert spaces,} London Mathematical Society Lecture
Note Series, 293, Cambridge University Press, Cambridge, 2002.
\bibitem{[H-Pe91]} Hu, Y., Peng, S. G., Adapted solution of a backward
semilinear stochastic evolution equation, \emph{Stochastic Anal.
Appl.,} 9, 4 (1991), 445--459.
\bibitem{[H-Pe96]} Hu, Y., Peng, S. G., Maximum principle for
optimal control of stochastic system of functional type,
\emph{Stochastic Anal. Appl.,} 14, 3 (1996),
283--301.
\bibitem{[Ichi]} Ichikawa, A., Stability of semilinear stochastic
evolution equations, \emph{J. Math. Anal. Appl.} 90, 1 (1982), 12--44.
\bibitem{[Lions]} Lions, J., \emph{Optimal control of systems governed by
differential equations,} Springer-Verlag, New-York, 1971.
\bibitem{[Oks05]} {\O}ksendal, B., Optimal Control of Stochastic Partial
Differential Equations, \emph{Stochastic Analysis and Applications,}
23 (2005), 165--179.
\bibitem{[Oks-Zh05]} {\O}ksendal, B., Proske, F., Zhang, T.,
Backward stochastic partial differential equations with jumps and
application to optimal control of random jump fields,
\emph{Stochastics,} 77, 5 (2005), 381--399.
\bibitem{[Pe-93]} Peng, S. G., Backward stochastic differential
equations and applications to optimal control, \emph{Appl. Math.
Optim.,} 27, 2 (1993), 125--144.
\bibitem{[Tess96]} Tessitore, G., Existence, uniqueness and space regularity
of the adapted solutions of a backward spde,
\emph{Stochastic Analysis and Applications,} 14, 4 (1996), 461--486.
\bibitem{[Y-Z]} Yong, J.,  Zhou, X. Y., \emph{Stochastic controls.
Hamiltonian systems and HJB equations,} Springer-Verlag, New-York,
1999.
\end{thebibliography}
}
\end{document}